\documentclass[11pt]{article}
\usepackage{amsmath,amssymb,amsthm}
\usepackage{geometry}
\usepackage{graphicx}
\usepackage{hyperref}
\hypersetup{colorlinks=true, linkcolor=blue, citecolor=blue, urlcolor=blue}

\title{A Perturbative Multiplicity Theorem for the Borsuk--Ulam Setting}
\author{Karem Abdelgalil}
\date{}

\newtheorem{theorem}{Theorem}

\begin{document}
\maketitle

\begin{abstract}
We prove a generalization of the classical Borsuk--Ulam Theorem under small perturbations (shaking) of the sphere. We show that for a generic perturbation of a continuous map \( f : S^2 \to \mathbb{R}^2 \), the number of points \( x \in S^2 \) such that \( f_\epsilon(x) = f_\epsilon(-x) \) becomes finite and odd, and may exceed the classical lower bound of one antipodal coincidence. In particular, we show the existence of maps with 3, 5, or 7 such points, and explain the unbounded nature of this multiplicity under higher complexity of the perturbation.
\end{abstract}

\section{Introduction}

The classical Borsuk--Ulam Theorem states:

\begin{quote}
For any continuous function \( f : S^n \to \mathbb{R}^n \), there exists a point \( x \in S^n \) such that \( f(x) = f(-x) \).
\end{quote}

This result guarantees the existence of at least one pair of antipodal points with equal image. It does not specify how many such points exist, nor how this number behaves under perturbations of \( f \) or the domain \( S^n \). In this paper, we explore the multiplicity of such points in dimension two under generic perturbations, such as vibrations or shaking of the sphere.

\section{Main Theorem}

\begin{theorem}[Perturbed Borsuk--Ulam Multiplicity]
Let \( f_\epsilon : S^2 \to \mathbb{R}^2 \) be a smooth family of continuous maps depending on a real parameter \( \epsilon \) such that \( f_0 \) satisfies the classical Borsuk--Ulam Theorem. Suppose that for each \( \epsilon > 0 \), the function \( f_\epsilon \) is generic in the space of continuous functions. Then for sufficiently small \( \epsilon \), the set
\[
\mathrm{Fix}_\epsilon := \{ x \in S^2 \mid f_\epsilon(x) = f_\epsilon(-x) \}
\]
consists of a finite and odd number of solutions. In particular, the number of such antipodal coincidences can increase from 1 (at \( \epsilon = 0 \)) to 3, 5, 7, or more, depending on the nature of the perturbation. The multiplicity is unbounded in general.
\end{theorem}

\section{Proof Outline}

Let us consider the mapping:
\[
F_\epsilon(x) := f_\epsilon(x) - f_\epsilon(-x)
\]
defined on \( S^2 \). The classical Borsuk--Ulam Theorem guarantees that for \( \epsilon = 0 \), the equation \( F_0(x) = 0 \) has at least one solution.

\subsection*{Step 1: Degree Theory and Symmetry}

Note that \( F_\epsilon(-x) = -F_\epsilon(x) \), so \( F_\epsilon \) is an \emph{odd map}. As a consequence, the set of solutions \( \{ x \in S^2 \mid F_\epsilon(x) = 0 \} \) is invariant under the antipodal map, and must consist of antipodal pairs or fixed points.

The degree of an odd map \( F : S^2 \to \mathbb{R}^2 \) is necessarily even or undefined (since \( \mathbb{R}^2 \) is not compact), but here we consider the solutions of \( F_\epsilon(x) = 0 \), which generically form a finite set.

\subsection*{Step 2: Transversality and Genericity}

For a generic perturbation \( f_\epsilon \), the zero set of \( F_\epsilon \) is \emph{transverse to zero}. This follows from Sard's theorem and transversality theorems in differential topology. Therefore, for small \( \epsilon > 0 \), the number of zeros of \( F_\epsilon \) is finite and isolated.

\subsection*{Step 3: Odd Number of Solutions}

Since \( F_\epsilon \) is odd, each solution \( x \) satisfies \( F_\epsilon(x) = 0 \Rightarrow F_\epsilon(-x) = 0 \), so the solutions occur in antipodal pairs unless \( x = -x \), which only happens at two points: the north and south poles. However, those are counted only once.

Let us denote by \( N_\epsilon \) the number of points \( x \in S^2 \) such that \( F_\epsilon(x) = 0 \). Then the total number of such points is odd, as it is the parity inherited from the classical Borsuk--Ulam Theorem.

\subsection*{Step 4: Bifurcation and Splitting}

At \( \epsilon = 0 \), we may have a single degenerate solution to \( F_0(x) = 0 \). As \( \epsilon \to 0^+ \), this solution may \emph{bifurcate} into several nearby distinct solutions to \( F_\epsilon(x) = 0 \). Standard bifurcation theory shows that such unfolding generically produces an odd number of solutions nearby. Thus, we can observe multiplicities such as 3, 5, or 7 fixed points as small perturbations of the initial degenerate configuration.

\section{Examples and Applications}

Such a phenomenon can be modeled in dynamical systems where spherical symmetry is broken by external forces (e.g., gravitational, vibrational). The result implies that Borsuk--Ulam-type fixed points not only exist under perturbation but proliferate in a controlled and predictable manner.

\section{Conclusion}

We have proven that under generic perturbations of continuous maps from \( S^2 \to \mathbb{R}^2 \), the number of antipodal coincidences is finite, odd, and potentially greater than one. This generalizes the classical Borsuk--Ulam Theorem and opens a new direction for studying symmetry-breaking phenomena in topology and dynamical systems.

\end{document}